\documentclass{article}

\usepackage[cp1251]{inputenc}
\usepackage[english]{babel}

\usepackage{amssymb,amsthm,amsmath}

\usepackage{graphicx}

\usepackage[all]{xy}

\renewcommand{\leq}{\leqslant}
\renewcommand{\geq}{\geqslant}

\begin{document}

\renewcommand{\figurename}{Pict.}

\begin{center}
{\large Makar Plakhotnyk}\\

{\Large \emph{On the history of the Ulam's Conjugacy}}\\

{\large \emph{(article for the journal for scholars)}}\\

{\large Postdoctoral researcher at the S\~ao Paulo University,
Brazil.\\ Mail: makar.plakhotnyk@gmail.com}
\end{center}

\begin{abstract}
We show the results on the history of the invention of the
conjugacy $h(x)=\frac{2}{\pi}\arcsin\sqrt{x}$ of one-dimensional
$[0,\, 1]\rightarrow [0,\, 1]$ maps $f(x)=4x(1-x)$ and
$g(x)=1-|1-2x|$.
\end{abstract}

\section{Introduction} Not every classical mathematical problem
has its final form in the moment of the creation. Moreover, this
is true about ``classical methods'' and ``classical tools'' of
solving mathematical problems. Almost all branches of modern
mathematics appeared or obtained the form, which is closed to
modern, in 18-19 century. During the reading of mathematical works
of that time, it is easy to regard such mistakes of bright
mathematicians, which should not be made by scholars and students
of firs cors of university nowadays. From another hand, it is
necessary to remember, that it is due to these (and only these)
mistakes the modern mathematics has the level of strictness, that
it has today.

We will pay attention to the appearance and the first attempts to
study the following construction. Let $f:\, I\rightarrow I$ be a
function, where $I$ is a set of real numbers. For any $n\in
\mathbb{N}$ denote $g(x) = \underbrace{f(f(\ldots f}_{n \text{
times}}(x)\ldots ))$, write $g = f^n$ for being short and call $g$
the $n$th iteration of $f$. At the beginning of the 19th century
the problems of finding the formula for $f^n$ by given $f$ or,
conversely, finding $f$ by given $g=f^n$ was formulated and
studied as deeply as it was possible. Non-succeed attempts to
solve these problems, lead to the appearance and active
development of Dynamical Systems Theory in the middle and second
part of the 20th century.

This article is motivated by the following question. It is known
that equality $g(x)=h(f(h^{-1}(x)))$ holds for every $x\in [0,\,
1]$, where $f(x)=4x(1-x)$, $g(x)=1-|1-2x|$ and
$h(x)=\frac{2}{\pi}\arcsin\sqrt{x}$ are defined on $[0,\, 1]$. The
mentioned equality appears in all textbooks on one-dimensional
dynamics in the chapter about topological conjugation (topological
equivalence). A lot of books on the Theory of Dynamical Systems
(for instance~\cite[p.~15]{Collet}, \cite[p.~103]{Hale-Kocak},
\cite[p.~316]{Katok} and~\cite[p.~341]{Melo-Strien}) refer this
formula to~\cite{Ulam}, but this is just a note on Summer Meeting
of the AMS in 1947, which is about one fifth of a page and
contains only the remark that $g$ can be used as a generator a
random numbers. During the private talks with Kyiv specialists on
one-dimensional dynamics the author failed to get the clearness
about the history of the investigation of this formula. All the
answers of professors looked like: "This is the mathematical
folklore: the fact is known to everybody, but not the authorship.
Certainly, one cay spend a lot of time and found some history, but
nobody cam be sure, wether or not this history will appear to be
interesting enough".

During the attempts to understand the history of the formula,
which is mentioned above, the author entered the magic world of
mathematical texts of 19-20 century, where the solutions of
mathematical problems are strongly mixed with ambitions of great
mathematicians and appear also at the pages of scientific
articles. But, let us look at everything in turn.

\section{Functional equations of J. Herschel}

John Herschel is a talented scientist of the end of 19th century.
For instance, he was presented with the Gold Medal of the Royal
Astronomical Society in 1826 and with the Lalande Medal of the
French Academy of Sciences in 1825, was a Honorary Member of the
St Petersburg Academy of Sciences. Also he was one of the founders
of the Royal Astronomical Society in 1820. Nevertheless, John
Herschel was an astronomer, but not a ``professional
mathematician''.

In the same time, J. Herschel stated in~\cite{Herschel} series of
mathematical problems and make attempts to solve them. His
solutions contain some disadvantages, which do not decrease the
value of the problems themselves.

\textbf{Problem 1.} \emph{Find $f^n$ for the function
$f(x)=2x^2-1$}. The following solution is presented. For every $x$
denote $u_0=x$ and consider the sequence
\begin{equation}\label{eq:01}
u_{n+1}=2u_n^2-1.\end{equation} To find the formula for $f^n$ is
the same as to find the general formula for $u_n$, dependent on
$n$ and $u_0$. Without any explanations Herschel decides to find
the solution in the form
\begin{equation}\label{eq:02} u_n =
\frac{1}{2}\left(C^{2^n}+C^{-2^n}\right),
\end{equation} where $C$ is a constant, dependent on $u_0$.
We may assume that the motivation for such form of the solution is
the formula $u_n=\lambda^n$ for the solution of a linear second
order difference equation $u_{n+2}=au_{n+1}+bu_n$. In any way, the
direct substitution lets to check that~\eqref{eq:02}
satisfies~\eqref{eq:01}. After the plug $u_0=x$
into~\eqref{eq:02}, Herschel has found $$ C=x+\sqrt{x^2-1},
$$ whence rewrites~\eqref{eq:02} as \begin{equation}\label{eq:06}
f^n(x)=\frac{1}{2}\left(\left(x+\sqrt{x^2-1}\right)^{2^n}+
\left(x-\sqrt{x^2-1}\right)^{2^n}\right).
\end{equation} He does not explain the possibility of
appearance of negative number under the square root.

\textbf{Problem 2.} \emph{Find functions $\varphi$, which satisfy
the equality}
\begin{equation}\label{eq:05}\varphi^2(x)=x.\end{equation} Like above,
J. Herschel denotes $x=u_n$ and $\varphi(u_n)=u_{n+1}$. Moreover,
he considers the sequence $\varphi(u_n)$, which denotes
$(\varphi(u))_n$. After the ``cross multiplication'' of equalities
\begin{equation}\label{eq:03} \left\{
\begin{array}{l}
(\varphi(u))_n=u_{n+1}\\
(\varphi(u))_{n+1}=u_n,
\end{array}\right.
\end{equation} he obtains $u_{n+1}(\varphi(u))_{n+1}
=u_n(\varphi(u))_n$, i.e. $u_n(\varphi(u))_n=c$ for all $n$, where
$c$ is a constant.

Also it follows from equalities~\eqref{eq:03} that
$u_n+(\varphi(u))_n= u_{n+1}+(\varphi(u))_{n+1}$, i.e.
$u_n+(\varphi(u))_n + C=0$ for all $n$ and some constant $C$. Now
Herschel claims that one may state $f(u_n(\varphi(u))_n)$ instead
of $C$, where $f$ is an arbitrary function. Thus the equation
\begin{equation}\label{eq:04} x + \varphi(x) +f(x\varphi(x))=0.
\end{equation} appears. Taking as example the function $f(x)=a+bx$
Herschel finds $\varphi$ from equality $$ x + \varphi(x) +
bx\varphi(x) +a=0
$$ as \begin{equation}\label{eq:16}
\varphi(x) = -\frac{a+x}{1+bx}
\end{equation}
and claims, that it is a solution of~\eqref{eq:05}. He does not
notice, that if one plug~\eqref{eq:16} into~\eqref{eq:05}, then
would not get the identity.

\textbf{Problem 3}. \emph{for given function $f$ find a function
$\varphi$ such that} $\varphi^n = f$. For this problem Herschel
suggests ``very simple'' solution: find the general formula for
$f^n$ and plug there $1/n$ instead of $n$. For example for the
function $f(x)=2x^2-1$ Herschel uses~\eqref{eq:06} to find
$$\varphi(x)=\frac{1}{2}\left(\left(x+\sqrt{x^2-1}\right)^{\sqrt[n]{2}}+
\left(x-\sqrt{x^2-1}\right)^{\sqrt[n]{2}}\right)$$ and remarks
that $\sqrt[n]{2}$ means the set of $n$th complex roots of $n$th
degree.

\textbf{Problem 4}. \emph{Let an hyperbola $AM$ with axis $CP$
(which coincides with $x$-axis) and center $C$ (which is origin)
be given. It is necessary to find the curve $am$ with the
following properties. For any point $P$ on the $x$-axis find the
point on $am$ with the same $x$-coordinate and take on the
$x$-axis the point $d$ such that $Cd = Pm$. Again find on $am$ the
point $l$ with $x$-coordinate as of $d$ and take the third point
$l_3$ such that $Cl_3 = dl$. Repeat these actions $n$ times, where
$n$ is chosen at the very beginning. The question: find $am$ such
that after $n$ steps we would get the segment $fk$, whose length
equals $PM$ (see fig.~\ref{fig:2}a.)}?

Write the equation of hyperbola $AM$ as $$ y^2 = (1-e^2)(a^2-x^2),
$$
where $e$ is the eccentricity of the hyperbola. For the function
$f(x)=\sqrt{(1-e^2)(a^2-x^2)}$ we have that
\begin{equation}\label{eq:07}
f^n(x)=\sqrt{(e^2-1)^nx^2
-\frac{e^2-1}{e^2-2}\left((e^2-1)^n-1\right)a^2}
\end{equation}
Moreover, finding a function $\varphi$ such that $\varphi^n=f$,
Herschel, according to his rule above, plugs $1/n$ instead of $n$
into~\eqref{eq:07} and obtains the equations of curves, which
satisfy the condition of the problem. He does not consider the
existence of another curves, which satisfy the former condition.
It is interesting, that Herschel does not notice here, that square
root has two values, one of which is complex.

\begin{figure}[htbp]
\begin{minipage}[h]{0.49\linewidth}
\begin{center}
\begin{picture}(140,100)

\put(10,10){\line(1,0){130}}

\qbezier(40,10)(45,54)(140,100)

\qbezier(60,10)(65,40)(140,80)

\qbezier(90,10)(90,10)(90,71.5)

\qbezier[25](90,71.5)(107.5,71.5)(125,71.5)

\qbezier(125,10)(125,10)(125,92)

\qbezier(107,10)(107,10)(107,83)

\put(43,12){$a$} \put(63,12){$A$} \put(12,12){$C$}

\put(88,0){$f$} \put(105,0){$d$}\put(123,0){$P$}

\put(88,74){$k$} \put(105,84){$l$} \put(123,96){$m$}
\put(127,66){$M$}

\end{picture}
\end{center}
\centerline{a. The picture from J. Herschel's work}
\end{minipage}
\hfill
\begin{minipage}[h]{0.49\linewidth}
\begin{center}
\begin{picture}(140,100)

\put(10,10){\line(1,0){130}}

\qbezier(30,10)(45,54)(140,100)

\put(28,0){$A$}

\qbezier(50,10)(50,10)(50,40)

\put(48,0){$B$} \put(48,44){$P$}

\qbezier(90,10)(90,10)(90,72.5)

\put(88,0){$C$}  \put(88,77){$Q$}

\qbezier(125,10)(125,10)(125,92.5)

\put(123,0){$F$} \put(123,95){$T$}

\end{picture}
\end{center}
\centerline{b. The picture from C. Babbage's work}
\end{minipage}
\caption{Graphical interpretations}\label{fig:2}
\end{figure}

\section{Conjugacy of  C. Babbage and J. Ritt}

Charles Babbage pays attention in~\cite[problem 9]{Babbage} to the
study of iterations of functions, precisely to the finding of the
function $\psi$ such that
\begin{equation}\label{eq:08} \psi^n(x)=x\end{equation} for all
$x$. First, he suggests the graphical interpretation of this
problem, which is similar to Herschel's one, which is described
above. At the picture~\ref{fig:2}b. it is necessary to find the
cure $APQT$ (where $A$ is the origin ) with the following
properties. For arbitrary point $B$ find the length of
perpendicular $BP$ and take on the $x$-axis the point $C$ such
that $AC = BP$. Then find the length of the perpendicular $CQ$ and
take its length on the $x$-axis. After $n$ steps we have to get
the perpendicular $FT$, whose length should be equal to the length
of the former segment $AB$. This geometrical construction is
exactly the verbal description of the equation~\eqref{eq:08}.

Babbage makes one more important remark. Let $f$ be some solution
of~\eqref{eq:08} and let $\varphi$ be an arbitrary invertible
function. Then
\begin{equation}\label{eq:09}g(x) =
\varphi(f(\varphi^{-1}(x)))\end{equation} will also be a solution
of~\eqref{eq:08}. In fact, these reasonings about~\eqref{eq:08} is
exactly the invention of the notion of topological conjugacy of
maps. Topological conjugacy of function $f$ and $g$ can be
imagined as the rectangle,
$$ 
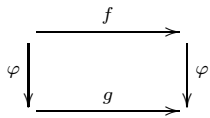
\xymatrix{
 \ar^{f}[rr] \ar_{\varphi}[d] &&
\ar^{\varphi}[d]\\
 \ar^{g}[rr] && }
$$ where we know, that one can come from the left top vertex to
the right bottom one either by ``top root'', or by ``bottom root''
with the same result. Coming by the root means here the sequent
applying the functions near the arrows to an arbitrary former
argument. These diagrams are called ``commutative diagrams'' and
are widely used in different branches of mathematics for the
illustration of reasonings. The important corollary of the
commutative diagram is the possibility ``to continue it to the
right'' and thus obtain
$$ 
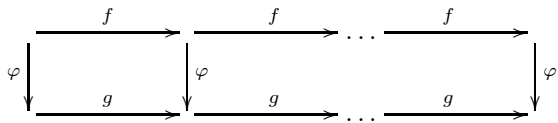
\xymatrix{
 \ar^{f}[rr] \ar_{\varphi}[d] &&
\ar^{\varphi}[d] \ar^{f}[rr] && \ldots \ar^{f}[rr] &&
\ar^{\varphi}[d]\\
 \ar^{g}[rr] &&  \ar^{g}[rr] && \ldots \ar^{g}[rr] && }
$$ which implies that $g^n(x) = \varphi(f^n(\varphi^{-1}(x)))$.
From another hand, the diagram just illustrates the reasonings and
the same conclusion could be obtained from the
equality~\eqref{eq:09} without the commutative diagram too.

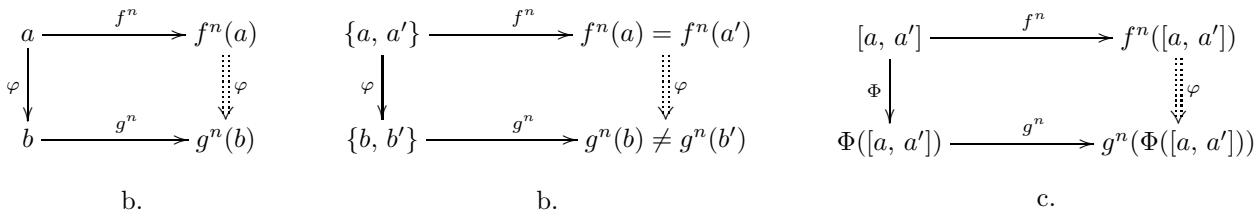
\begin{figure}[htbp]
\begin{minipage}[h]{0.22\linewidth}
$$ \xymatrix{
a \ar^{f^n}[rr] \ar_{\varphi}[d] &&
f^n(a) \ar@{:>}^{\varphi}[d]\\
b \ar^{g^n}[rr] && g^n(b)}
$$
\centerline{b.}
\end{minipage}
\hfill
\begin{minipage}[h]{0.35\linewidth}
$$ \xymatrix{
\{a,\, a'\} \ar^{f^n}[rr] \ar_{\varphi}[d] &&
f^n(a)=f^n(a') \ar@{:>}^{\varphi}[d]\\
\{b,\, b'\}  \ar^{g^n}[rr] && g^n(b)\neq g^n(b')}
$$
\centerline{b.}
\end{minipage}
\hfill
\begin{minipage}[h]{0.35\linewidth}
$$ \xymatrix{
[a,\, a'] \ar^{f^n}[rr] \ar_{\Phi}[d] &&
f^n([a,\, a']) \ar@{:>}^{\varphi}[d]\\
\Phi([a,\, a'])  \ar^{g^n}[rr] && g^n(\Phi([a,\, a']))}
$$
\centerline{c.}
\end{minipage}
\caption{Reasonings of Ritt}\label{fig:3}
\end{figure}

Joseph Ritt (injustively) says in~\cite{Ritt} that C. Babbage
earlier claimed that for any fixed solution $f$ of the
equation~\eqref{eq:08} and any other solution $g$ there exists an
invertible $\varphi$ such that~\eqref{eq:09} holds. Suppose
(writes Ritt) that $\varphi(a)=b$ for some $a$ and $b$. Then
$\varphi(f^n(a))=g^n(b)$ for every $n\in \mathbb{N}$. If for some
other $a',\, b'$ the equality $\varphi(a')=b'$ holds, then for
every $n$ we will have $\varphi(f^n(a'))=g^n(b')$ (see.
pict.~\ref{fig:3}a). If $f^n(a)=f^n(a')$ for some $n$, but
$g^n(b)\neq g^n(b')$, the contradiction with that $\varphi$ is a
function would appear (see pict.~\ref{fig:3}b).

Ritt continues these reasonings. Suppose that $\varphi(x)=\Psi(x)$
for all $x\in (a,\, a')$ and some $a,\, a'$, where $\Psi$ is some
function. In other words, suppose that $\varphi$ is already
defined on some interval $(a,\, a')$. Then for any $n$ the
equality $\varphi(f^n(x))=g^n(\varphi((x)))$ determines $\varphi$
on $f^n([a,\, a'])$ (see pict.~\ref{fig:3}c). If for some $n_1,\,
n_2$ intervals $f^{n_1}([a,\, a'])$ and $f^{n_2}([a,\, a'])$
intersect, the the contradiction may appear.

Notice, that talking about the contradiction, Ritt does not pay
attention to the fact, which gives the contradiction. Indeed, he
assumes that two points $(a,\, b)$ and $(a',\, b')$ such that
$\varphi(a)=b$ and $\varphi(a')=b'$ are given and then
contradiction appears. But this only means that the graph of
$\varphi$ does not pass through both of these points (but may be
pass through one of them). It is not clear, what can generate the
conclusion that for some exact solutions $f$ and $g$
of~\eqref{eq:08} there in no any invertible $\varphi$, which
transforms~\eqref{eq:09} to the identity.

\section{Lamerey's Diagrams}

The Cartesian method as a way on plotting the graphs of functions
is known from the de Cartes time, i.e. from the first part of the
17th century. In the same time, the problems, which were stated by
J. Herschel (and were discussed above), need a bit specifical
techniques, and we will pay attention to it now.

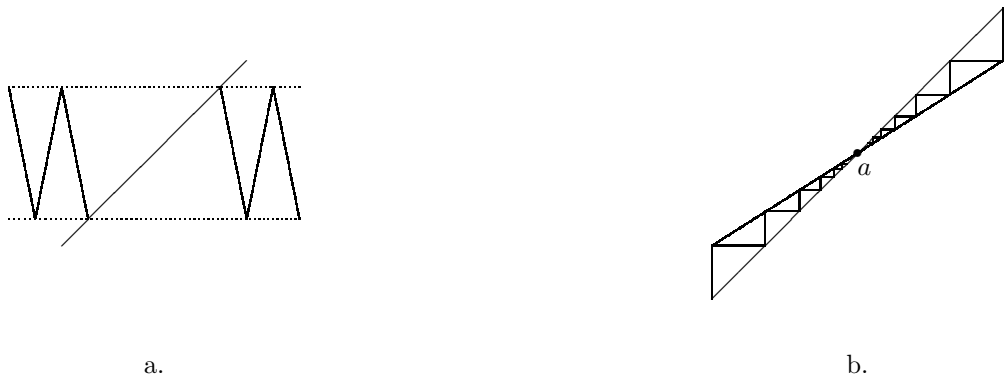
\begin{figure}[htbp]
\begin{minipage}[h]{0.45\linewidth}
\begin{center}
\begin{picture}(110,110)

\put(20,20){\line(1,1){70}}

\qbezier[60](0,30)(60,30)(110,30)
\qbezier[60](0,80)(55,80)(110,80)

\qbezier(30,30)(25,55)(20,80) \qbezier(10,30)(15,55)(20,80)
\qbezier(0,80)(5,55)(10,30)

\qbezier(80,80)(85,55)(90,30) \qbezier(90,30)(95,55)(100,80)
\qbezier(100,80)(105,55)(110,30)

\end{picture}
\end{center}
\centerline{a.}
\end{minipage}
\hfill
\begin{minipage}[h]{0.45\linewidth}
\begin{center}
\begin{picture}(110,110)
\put(0,0){\line(1,1){110}}

\qbezier(0,20)(55,55)(110,90)

\put(55,55){\circle*{3}} \put(55,47){$a$}

\qbezier(0,0)(0,10)(0,20) \qbezier(0,20)(10,20)(20,20)
\qbezier(20,20)(20,26)(20,33) \qbezier(20,33)(26,33)(33,33)
\qbezier(33,33)(33,37)(33,41) \qbezier(33,41)(37,41)(41,41)
\qbezier(41,41)(41,43)(41,46) \qbezier(41,46)(43,46)(46,46)
\qbezier(46,46)(46,47)(46,49) \qbezier(46,49)(47,49)(49,49)
\qbezier(49,49)(49,50)(49,51) \qbezier(49,51)(50,51)(51,51)

\qbezier(59,59)(60,59)(61,59) \qbezier(61,59)(61,60)(61,61)
\qbezier(61,61)(62,61)(64,61) \qbezier(64,61)(64,62)(64,64)
\qbezier(64,64)(66,64)(69,64) \qbezier(69,64)(69,66)(69,69)
\qbezier(69,69)(73,69)(77,69) \qbezier(77,69)(77,73)(77,77)
\qbezier(77,77)(83,77)(90,77) \qbezier(90,77)(90,83)(90,90)
\qbezier(90,90)(100,90)(110,90) \qbezier(110,90)(110,105)(110,110)
\end{picture}
\end{center}
\centerline{b.}
\end{minipage}
\caption{Lamerey's Diagrams}\label{fig:4}
\end{figure}

When one use the Cartesian method, the arguments of a functions
are understood as points of the $x$-axis and its values -- as
points of $y$-axis, whence the function appears to be a map from
$x$-axis to $y$-axis. In the same time, we can consider any
function $f$ as a map from the line $y=x$ to itself, which maps a
point $(x,\, x)$ to $(f(x),\, f(x))$. The graphical interpretation
of this action can be imagined as vertical line from point $(x,\,
x)$ to the graph and then horizontal line to $(f(x),\, f(x))$,
returning to $y=x$.

For usefulness of such interpretation, let us come back to the
problem of finding a continuous function $\varphi:\,
\mathbb{R}\rightarrow \mathbb{R}$ such that
$\varphi^2(x)=\varphi(x)$ for all $x$ from the domain. Clearly,
for every $x_0$ the equality $\varphi(\varphi(x_0))=\varphi(x_0)$
holds. Thus, denote $y_0=\varphi(x_0)$ and obtain that
$\varphi(y_0)=y_0$, i.e. the graph of $\varphi$ passes through the
point $(y,\, y)$ for every $y$ such that $y=\varphi(x)$ for some
$x$. This means that there exist $a,\, b$ such that
$\varphi(\mathbb{R}) = [a,\, b]$ and $\varphi(x)=x$ for every
$x\in [a,\, b]$ (see pict.~\ref{fig:4}a).

Also the following reasonings can be illustrated. Notice that the
solutions of the equation $x=\varphi(x)$ are precisely the points
of the intersection of the graph of $\varphi$ and the line $y=x$.
Suppose that for the function $\varphi$ there exists (unknown)
solution $a$ of the equation $x=\varphi(x)$ and, moreover,
$x<\varphi(x)<a$ for every $x<a$ and $a<\varphi(x)<x$ for each
$x>a$. For arbitrary $x_0$ consider the sequence $x_k =
\varphi^k(x_0)$. If one plot the graph, then it becomes evident
that the sequence $\{x_k\}$ tends to unknown $a$ (see
pict.~\ref{fig:4}b). This reasonings were invented in the
beginning of the 19th century almost simultaneously by Evarist
Galois~\cite{Galois} and Adrien Legendre~\cite{Legendre} (see
also~\cite{Days}).

\begin{figure}
\begin{minipage}[h]{0.45\linewidth}
\begin{center}
\begin{picture}(140,140)

\put(10,10){\line(1,0){130}} \put(10,10){\line(0,1){130}}
\put(10,10){\line(1,1){120}} \put(10,130){\line(1,0){120}}
\put(130,10){\line(0,1){120}} \put(10,10){\line(0,1){120}}

\linethickness{0.3mm} \qbezier(10,10)(10,10)(34,58)
\qbezier(34,58)(34,58)(82,106) \qbezier(82,106)(82,106)(130,130)

\qbezier[20](10,10)(16,44)(22,58)
\qbezier[15](22,58)(28,70)(34,82)
\qbezier[20](34,82)(46,94)(58,106)
\qbezier[15](58,106)(70,112)(82,118)
\qbezier[23](82,118)(106,124)(130,130) \linethickness{0.1mm}

\qbezier(22,10)(22,10)(22,58) \qbezier(22,34)(22,34)(34,34)

\qbezier(34,10)(34,10)(34,82) \qbezier(34,82)(34,82)(82,82)

\qbezier(22,58)(22,58)(58,58) \qbezier(58,10)(58,10)(58,106)
\qbezier(58,106)(58,106)(106,106) \qbezier(82,10)(82,10)(82,118)
\qbezier(106,106)(106,106)(106,118)
\qbezier(82,118)(82,118)(106,118)

\put(20,2){$x_1$} \put(32,2){$x_2$} \put(56,2){$x_3$}
\put(80,2){$x_4$}

\end{picture}
\end{center}
\end{minipage}
\hfill
\begin{minipage}[h]{0.45\linewidth}
\begin{center}
\begin{picture}(140,140)

\put(10,10){\line(1,0){130}} \put(10,10){\line(0,1){130}}
\put(10,10){\line(1,1){120}} \put(10,130){\line(1,0){120}}
\put(130,10){\line(0,1){120}} \put(10,10){\line(0,1){120}}

\qbezier(30,54)(30,111)(30,111) \qbezier(30,54)(30,54)(54,54)
\qbezier(54,54)(54,54)(54,111) \qbezier(54,111)(54,111)(30,111)

\qbezier(70,10)(70,120)(70,119.5)

\linethickness{0.3mm}

\qbezier(10,10)(10,130)(130,130) \qbezier(10,10)(10,60)(70,70)
\qbezier(70,70)(120,80)(130,130)

\qbezier[40](10,10)(11,80)(30,111)
\qbezier[25](30,111)(60,119)(70,119.5)
\qbezier[40](70,119.5)(105,120)(130,130) \linethickness{0.1mm}

\put(70,119.5){\circle*{4}} \put(30,111){\circle*{4}}

\put(100,73){$f$} \put(40,65){$f$}

\put(40,93){$\varphi$}

\put(35,120){$\varphi(f)$}

\put(30,45){$1$} \put(55,45){$2$} \put(56,103){$3$}
\put(22,113){$4$}

\end{picture}
\end{center}
\end{minipage}
\caption{Pictures from the S. Pincherle's work} \label{fig:1}
\end{figure}
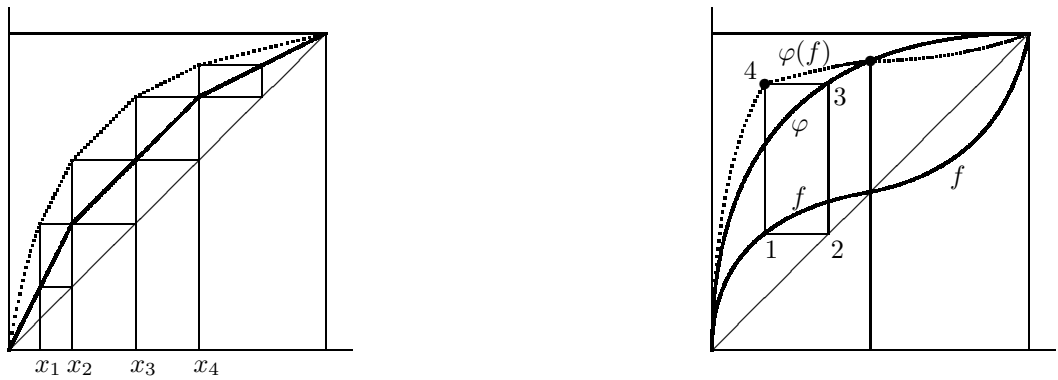

The work~\cite{Pincherle} contains the picture (see
pict.~\ref{fig:1}), where this method is illustrated. The left
graph contains a line $y=x$ (for the technical calculations,
mentioned earlier) and the graph $y=f(x)$, which is used to
construct the function $y=f^2(x)$. The right sketch contains the
construction of the graph of the function $y(x)=\varphi(f(x))$ by
given graphs of $f$ and $\varphi$.

Nowadays the sketches like those on figure~\ref{fig:1}, are called
Lamerey's Diagrams due to the work~\cite{Lemerai}.

\section{Trigonometry of J. Boole}

Gorge Boole is known in the history of mathematics as one on the
foundators of the mathematical logics. In~\cite{Boole} G.~Boole
suggests the following way of finding the general formula for
iterations of the function $f(x)=2x^2-1$. Keeping in mind the
double angle formula $\cos 2x = 2\cos^2x-1$, denote $g(x)=2x$ and
$h(x)=\cos x$, whence write the commutative diagram
\begin{equation}\label{eq:14} \xymatrix{ x \ar^{g}[rr] \ar_{h}[d]
&& 2x
\ar^{h}[d]\\
\cos x \ar^{f}[rr] && f(h)=h(g)  } \end{equation} Continue the
diagram to the right and obtain
$$
\xymatrix{ x \ar^{g}[rr] \ar_{h}[d] && 2x\ar^{h}[d] \ar^{g}[rr] &&
\ldots \ar^{g}[rr] && 2^nx
\ar^{h}[d]\\
\cos(x) \ar^{f}[rr] && \cos(2x) \ar^{f}[rr] && \ldots \ar^{f}[rr]
&& \cos(2^nx),}
$$
whence $f^n(\cos x) = \cos(2^nx).$ The substitution $t=\cos x$
leads to
$$f^n(t) = \cos(2^n\arccos t)$$ for all $t\in [-1,\, 1]$.

\section{The J.
von Neumann's and S. Ulam's generator}

John von Neumann studied deeply quantum physics,functional
analysis, sets theory and informatics. His name is connected with
the architecture of the most modern computers. von Neumann,
together with polish mathematician Stanislaw Ulam participated in
Manhattan Project.

In the short note~\cite{Ulam}, which is the thesis of the mating
of American Mathematical Society, J. von Neumann and S. Ulam
suggested to use iterations of the function
\begin{equation}\label{eq:13} f(x)=4x(1-x)\end{equation}
for obtaining ``numbers with different distributions''.
Nevertheless, they have not explained there, what the
``distribution of a number'' is.

The complicatedness of calculations, which are dealing with the
function~\eqref{eq:13}, is known from the Pier
Verhulst's~\cite{Verhulst-1}, written in the first part of the
19th century, where he used the formula
\begin{equation}\label{eq:15}
p_{k+1} = p_k(m-np_k).
\end{equation} for the expectation of the number of individuals $p_k$ in
the biological population in the $k$th generation, where $m$ and
$n$ are constants. Verhulst remarked the strong dependence  of
$p_k$ on $p_1,\, m$ and $n$. In other words, for the huge $k$ the
small changes of $p_1,\, m$ and $n$ can lead to march more change
of $p_k$. The sequence~\eqref{eq:15} was called ``logistic'' due
to ``logists'', who made calculations in the Ancient Greece. The
map~\eqref{eq:13} is also called logistic map, because it is of
the form~\eqref{eq:15}.

In the article~\cite{Neumann} von Neumann mentioned the
disadvantage of the function~\eqref{eq:13} as a generator,
introduced in~\cite{Ulam} and, in the same time, explained, what
is meant under the distribution of a number. Suppose that
non-decreasing function $F:\, [0,\, 1]\rightarrow [0,\, 1]$ such
that $F(0)=0$ and $F(1)=1$ is given. It is necessary ``in some
way'' to get the sequence of numbers $\{x_n\}\subset [0,\, 1]$
such that the probability of the event $x_i\in [a,\ b]$ equals
$F(b)-F(a)$. The result of~\cite{Ulam}, says von Neumann, is that
for almost every $x_0\in [0,\, 1]$ (up to Lebesgue measure), the
probability $x_i\in [a,\, b]$ equals $b-a$. He mentions that
obtaining the function, which ``generates'' the sequences with
given distribution, is an important problem, for instance, for
some calculating methods. In the same time, it is mentioned
in~\cite{Neumann} that the function, suggested in~\cite{Ulam} can
not be used with the mentioned goal. The arguments for such
impossibility were the following.

For a given sequence $\{x_i\}\subset [0,\, 1]$ such that
$x_{i+1}=4x_i(1-x_i)$ define $\alpha_i$ by
$x_i=\sin^2\pi\alpha_i$. Now notice that $\alpha_{i+1}=2\alpha_i\,
(mod\ 1)$, i.e. $\alpha_{i+1}$ is the fractional part of
$2\alpha_i$. Denote the binary decomposition of $\alpha_1$ as
$\alpha_1=\beta_1\beta_2\beta_3\ldots $, whence the binary
decomposition of $\alpha_k$ would be
$\alpha_k=\beta_k\beta_{k+1}\beta_{k+2}\ldots $. Now von Neumann
makes conclusion that in real computer we can not take a number
with infinitely many random binary digits, whence the sequence
$\{\alpha_i\}$ (and, correspondingly, $\{x_i\}$), being considered
``on the real computer'' will become $0$ after finitely many
iterations.

Moreover, von Neumann made in~\cite{Neumann} a very non-ethical
think: he wrote that is was S. Ulam, who suggested~\eqref{eq:13}
as a random generator and, thus, ``forget'' that Ulam was
co-author of~\cite{Ulam} too. The continuation of this history a
much more interesting. Is it easy no see that the formula
$x_i=\sin^2\pi\alpha_i$ does not define the sequence $\alpha_i$
with properties, which are mentioned in~\cite{Neumann}. If one
would find $\alpha\in [0,\, 1]$ then there is no one-to-one
correspondence, because $\alpha_i = k+\frac{(-1)^k}{\pi}\arcsin
\sqrt{x_i},\, k\in \{0,\, 1\}$ and it is not clear, it is
necessary to take $k$ being $0$ or $1$. If suppose that
$\alpha_i\in [0,\, 0.5]$, then the formula
$\alpha_{i+1}=2\alpha_i\, (mod\ 1)$ would give a number from the
interval $(0.5,\, 1]$ after finitely many steps for every
$\alpha_1\neq 0$.

Remark, that the function $\alpha(x) =
\frac{1}{\pi}\arcsin{\sqrt{x}}$ is the bijection between the
intervals $[0,\, 1]$ and $[0,\, 0.5]$. Thus, commutative diagram
$$ \xymatrix{ [0,\, 1] \ar^{f}[rr] \ar_{\alpha}[d] && [0,\, 1]
\ar^{\alpha}[d]\\
[0,\, 0.5] \ar^{g}[rr] && [0,\, 0.5]  }
$$ defines a function $g:\, [0,\, 0.5]\rightarrow [0,\, 0.5]$ by
$g(x)=\alpha(f(\alpha^{-1}(x)))$. Remind that $\alpha^{-1}(x) =
\sin^2\pi x$. The evident technical calculations give the
following:
$$ g(x)=\frac{1}{\pi}\arcsin\sqrt{4\sin^2(\pi x)(1-\sin^2\pi x)}
=\frac{\arcsin(|\sin (2\pi x)|)}{\pi}.
$$ %
Simplification of the absolute value function in the argument of
$\arcsin$ leats to the dichotomy either $x\leq 0.25$ or $x>0.25$
(this is naturally, since we have already seen some ``problems''
with the middle of the interval for $x$ in the von Neumann's
work). In fact,
\begin{equation}\label{eq:10} g(x) = \left\{
\begin{array}{ll}
2x& \text{for  }x\in [0,\, 0.25]\\
1-2x& \text{for }x\in (0.25,\, 0.5]
\end{array}\right.
\end{equation}

We have mentioned already, that commutative diagrams are connected
with the conjugacy. Let us give the formal definition. Functions
$f:\, A\rightarrow A$ and $g:\, B\rightarrow B$, where $A$ and $B$
are sets of real numbers, are called topologically conjugated, if
there exists an invertible function $h:\, A\rightarrow B$ such
that $h(f(x))=g(h(x))$ for all $x\in A$.

The chapter of mathematics ``topology'' studies geometrical
``figures'' up to some transformations, which roughly can be
imagined like if figures be made from a material, which admits
stretches and compressing. The word combination ``topological
conjugation'' can be explained as the change of the scale while
the construction of the graph of a function. We shall explain more
carefully what ``the change of the scale'' is. Let us construct
the graph of the function $f(x)=4x(1-x)$ for $x\in [0,\, 1]$ (it
is the parabola inside the square $[0,\, 1]\times [0,\, 1]$, with
branches going down). Let $h:\, [0,\, 1]\rightarrow [a,\, b]$ be a
continuous invertible function, which (we will say defines the
change of coordinates). Let us take each point on the segment
$[0,\, 1]$ of $x$-axis and $y$-axis and write $h(x)$ near it.
After this, we may say, that each point $(x,\, y)$ of the square
$[0,\, 1]\times [0,\, 1]$ obtains the new coordinates
$(\widetilde{x},\, \widetilde{y})\in [a,\, b]\times [a,\, b]$, but
the line $y=x$ remains to be $y=x$. Nevertheless, the parabola,
which was (!) the graph of the function $f(x)=4x(1-x)$ appears to
the the graph of some another function $\widetilde{f}:\, [a,\,
b]\rightarrow [a,\, b]$. If $\alpha:\, [0,\, 1]\rightarrow [0,\,
0.5]$ would be instead of $h$ above, then it follows from our
calculation, that the obtained ``new function'' is $g$ of the
form~\eqref{eq:10}.

After $g$ has been found, it is necessary to ``want'' to obtain
one more function, which will be ``similar'' to $g$, but will be
$[0,\, 1]\rightarrow [0,\, 1]$ instead of $[0,\, 0.5]\rightarrow
[0,\, 0.5]$. If apply the topological conjugacy $\beta(x)=2x$ to
$g$, then obtain $g_1(x)=\beta(g(\beta^{-1}(x)))$, which is
\begin{equation}\label{eq:11} g_1(x) = \left\{
\begin{array}{ll}
2x& \text{for }x\in [0,\, 0.5],\\
2-2x& \text{for }x\in (0.5,\, 1].
\end{array}\right.
\end{equation}%
Notice that graphs of functions~\eqref{eq:10} and~\eqref{eq:11}
looks the same if the numbers are marked uniformly on coordinate
axis. The conjugation of $f$ and $g_1$ can be easily illustrated
by the following commutative diagram $$ \xymatrix{ [0,\, 1]
\ar^{f}[rr] \ar_{\alpha}[d] \ar@/_3pc/@{-->}_{h}[dd] && [0,\, 1]
\ar@/^3pc/@{-->}^{h}[dd]
\ar^{\alpha}[d]\\
[0,\, 0.5] \ar^{g}[rr] \ar_{\beta}[d] && [0,\, 0.5]
\ar^{\beta}[d]\\
[0,\, 1] \ar^{g_1}[rr] && [0,\, 1]  }
$$ %
and the conjugacy (the function, which defines the conjugation) is
\begin{equation}\label{eq:12} h(x)= \beta(\alpha(x)) =
\frac{2}{\pi}\, \arcsin\sqrt{x}. \end{equation}

Let us come bach to the George Boole's work~\cite{Boole}.
Commutative diagram~\eqref{eq:14} does not define a topological
conjugation, because the map $g$ there is not $B\rightarrow B$
such that $h(x)=\cos x$ is the bijection between $B$ and $h(B)$.

\begin{figure}[htbp]
\begin{minipage}[h]{0.45\linewidth}
\begin{center}
\begin{picture}(110,110)
\put(0,0){\line(1,1){110}} \put(0,0){\vector(1,0){110}}
\put(0,0){\vector(0,1){110}}

\qbezier(0,0)(22.5,45)(45,90) \qbezier(45,90)(67.5,45)(90,0)

\qbezier(0,0)(45,180)(90,0)

\qbezier[60](0,90)(45,90)(90,90) \qbezier[60](90,90)(90,45)(90,0)

\end{picture}
\end{center}
\centerline{a.}
\end{minipage}
\hfill
\begin{minipage}[h]{0.45\linewidth}
\begin{center}
\begin{picture}(110,110)
\put(0,0){\line(1,1){110}} \put(0,55){\vector(1,0){110}}
\put(55,0){\vector(0,1){110}}

\qbezier(10,100)(50,-80)(100,100)
\qbezier[60](10,100)(32.5,55)(55,10)
\qbezier(55,10)(77.5,55)(100,100)

\qbezier[60](10,100)(55,100)(100,100)
\qbezier[60](100,100)(100,55)(100,10)
\qbezier[60](10,10)(55,10)(100,10)
\qbezier[60](10,100)(10,55)(10,10)

\end{picture}
\end{center}
\centerline{b.}
\end{minipage}
\caption{Topologically conjugated maps}\label{fig:5}
\end{figure}
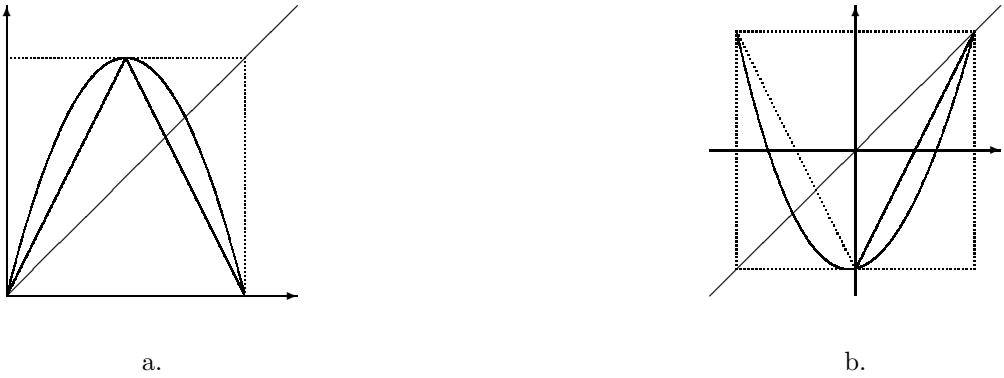

Graphs of~\eqref{eq:13} and~\eqref{eq:11} are given at
picture~\ref{fig:5}a. Picture~\ref{fig:5}b contains graphs of
$y=2x^2-1$ and $y=2x$, which appeared in the Boole's~\cite{Boole}.
It is clear, that these picture are ``geometrically'' equal, just
one of them is turned upside-down from another.

The conjugacy $h(x)=1-x$ turns the parabola from the
Picture~\ref{fig:5}a to the form of Picture~\ref{fig:5}b. The fact
that coordinate axis are ``at another position'' mean nothing,
because the conjugacy $h(x)=x+a$, where $a$ is a constant, moves
the graph along the line $y=x$ (the same line, whose importance we
have already mentioned, talking about the iterations). Thus, both
George Boole in~\cite{Boole} and John von Neumann
in~\cite{Neumann}, studied the iterations of a parabola (it is not
important, which parabola) and almost invented the topological
conjugation of a parabola and a piecewise linear function, but
stopped with the consideration of just one linear part of it.

Also Stanislaw Ulam returns in~\cite{Ulam-1964} to the question,
which we have discussed at the end of the chapter, dedicated to
Joseph Ritt's work. Let $g_1:\, [0,\, 1]\rightarrow [0,\, 1]$ be
of the form~\eqref{eq:11} and $g_2:\, [0,\, 1]\rightarrow [0,\,
1]$ be defined as follows:
$$g_2(x) = \left\{
\begin{array}{ll}
l(x)& \text{for }x\in [0,\, v],\\
r(x)& \text{for }x\in (v,\, 1],
\end{array}\right.
$$ where $v\in (0,\, 1)$ is a parameter, function $l$ increase, $r$
decrease such that $l(0)=r(1)=0$. In other words, the graph of
$g_2$ is ``similar'' to the graph of $g_1$, but $g_2$ is not
necessary piecewise linear and is not necessary symmetrical in
$x=0.5$. Ulam asks: which conditions should satisfy $g_2$ for
being topologically equivalent to $g_1$ (or, which is the same, to
$y=4x(1-x)$)? Suppose that $\widetilde{h}:\, [0,\, 1]\rightarrow
[0,\, 1]$ is such that
$g_2(x)=\widetilde{h}(g_1(\widetilde{h}^{-1}(x)))$ for all $x\in
[0,\, 1]$. Since $\widetilde{h}$ is invertible, then
$\widetilde{h}(0)\in \{ 0,\, 1\}$. Now $g_2(0)=0$ and
$\widetilde{h}(g_1(0))=g_2(\widetilde{h}(0))$ imply
$\widetilde{h}(0)=0$, whence $\widetilde{h}$ increase. Moreover,
since $g_1^n(0)=0$ for all $n\geq 1$, then it follows from
$\widetilde{h}(g_1^n(x))=g_2^n(\widetilde{h}(x))$ that
$g_2^n(\widetilde{h}(x))=0$, whenever $g_1^n(x)=0$. Ulam makes a
conclusion from these reasonings that $g_2$ is topologically
conjugated to $g_1$ is and only if that set $M = \{ x\in [0,\, 1]
:\, g_2^k(x)=0\ \text{ for some }k\}$ is dense in $[0,\, 1]$.

As about the map $\psi$, which satisfies the functional
equation~\eqref{eq:08}, Melvyn Nathanson has proved the following
theorem in~\cite{Natanson}: if $\psi:\, [0,\, 1]\rightarrow [0,\,
1]$ is a continuous function such that $\psi^p(x)=x$ for all $x\in
[0,\, 1]$, then $\psi^2(x)=x$ for all $x\in [0,\, 1]$. Moreover,
if $p$ is even, then $\psi(x)=x$ for all $x\in [0,\, 1]$. In other
words, the problem, which Babbage treated, has, in sone cense, the
only trivial solution. Nevertheless, his does not decrease the
importance or the reasonings, which were suggested by Babbage
during the attempts of the solution. The article~\cite{Natanson}
is written in the second part of the 20th century and is dedicated
to the Li and Yorke's  notion of chaos.

The formula~\eqref{eq:12} or the topological conjugacy of
maps~\eqref{eq:13} and~\eqref{eq:11} was published at first by
Ottis Richard in~\cite{Rechard}. Nevertheless, Richard thanks
there S. Ulam ``for many helpful and stimulating conversations on
the subject of this paper'' and writes especially about the
formula~\eqref{eq:12} (which is used for the calculating of the
invariant measure) that it is S. Ulam, which noticed this formula
at first. Stanislaw Ulam published~\eqref{eq:12} the first time 8
years later in~\cite{Ulam-1964}. The article~\cite{Ulam-1964} does
not contain any references to John von Neumann as the author of
idea of getting~\eqref{eq:12}. Some modern books on Dynamical
Systems theory refer~\eqref{eq:12} as Ulam's map.

\section{Final remarks}

Sometimes mathematicians make mistakes. Mathematicians of the
worldwide level sometimes make mistakes too. Sometimes
mathematicians quarrel one with each other and tracks of this can
be seen in their articles. Nevertheless, the mistake of a scholar
and the mistake of a mathematician of the worldwide level are
``different sings''.

John Herschel's problems with solutions, which were published in
year 1814, contain mistakes. But they also contain the foundations
of the theory of one-dimensional dynamical systems - a theory,
which was developed in the second part of 20th century. The idea,
which lead to the notion of topological equivalency, is clearly
formulated in year 1815 during the attempts of the description of
a class of functions, which, in some thence, is empty (but it is
the result of 1970th).

Lamerey's Diagrams, which were described in 1897, were used, in
fact, by Adrien Legendre in 1808 and Evariste Galois in 1830. John
von Neumann refuse his participation in the preparation of the
work, which contained a mistake, but this lead to that he loosed
the authority for the more important result, which is mentioned
now without the name of von Neumann.

Notice, that pictures~\ref{fig:2} (about Herschel's and Babbage's
works) and~\ref{fig:1} (about Pincherle's work) are prepared by
the graphical environment of \LaTeX, i.e. they are not copies of
the original articles, but saves the meaning and the notations of
the original ones.

\end{document}